\documentstyle[amssymb,11pt]{article}
\begin{document}

\begin{center}{\bf EQUATIONS OF MOTION WITH RESPECT TO THE $(1+1+3)$ THREADING OF A $5D$ UNIVERSE} \end{center}
\begin{center}{AUREL BEJANCU \\ Department of Mathematics \\ Kuwait University\\ P.O.Box 5969, Safat 13060\\ Kuwait\\ 
 E-mail:aurel.bejancu@ku.edu.kw}\end{center}

\begin{abstract}
{We continue our research work started in [1], and obtain in a covariant form, the equations of motion with respect
 to the $(1+1+3)$ threading of a $5D$ universe $(\bar{M}, \bar{g})$. The natural splitting of the tangent bundle 
of $\bar{M}$ leads us to the study of three categories of geodesics: spatial geodesics, temporal geodesics and 
vertical geodesics. As an application of the general theory, we introduce and study what we call the $5D$ 
 Robertson-Walker universe.}
 \end{abstract}

%\noindent {\bf Keywords}: (1+3) threading of spacetimes; kinematic quantities; Raychaudhuri's equation; 
%Ricci tensor field; Riemannian spatial connection; spatial tensor fields.

%PACS number(s): 04.50.-h
 
\medskip
\newpage
\section{Introduction}
This paper is a continuation of our previous paper \cite{b} on kinematic quantities and Raychaudhuri equations in 
 a $5D$ universe. According to the new approach presented in \cite{b}, the $5D$ universe $\bar{M} = M \times K$ 
is studied by means of the submersion of $\bar{M}$  on the $4D$ spacetime $M$.  Note that in all the other theories 
of a $5D$ universe the study was performed via an immersion of $M$ in $\bar{M}$ (cf.[2-5]).\par
 The kinematic quantities together with the spatial tensor fields and the Riemannian spatial connection enable us
 to obtain, in a covariant form, the equations of motion in ($\bar{M}, \bar{g}$). By using the natural splitting of
 the tangent bundle of $\bar{M}$ we introduce into the study three categories of geodesics: spatial geodesics, 
 temporal geodesics and vertical geodesics. We apply the general theory to what we call the $5D$ Robertson-Walker 
 universe, which can be thought as a disjoint union of $4D$ Robertson-Walker spacetimes. In this case, the above
 three categories of geodesics are completely determined.\par
Now, we outline the content of the paper. In Section 2 we recall from \cite{b} the kinematic quantities with
 respect to the $(1 + 1+ 3)$ threading of the $5D$ universe ($\bar{M}, \bar{g}$), and the Riemannian spatial 
connection $\nabla$ on the spatial distribution ${\cal{S}}\bar{M}$. The complete characterization of the 
Levi-Civita connection on $(\bar{M}, \bar{g})$  (cf. (2.18)) enables us to write down in Section 3,  for the first
 time in literature, the splitting of the equations of motions into three groups (cf.(3.6)). As an example, we
 present the $5D$ Robertson-Walker universe (see (3.7)) together with its equations of motion (cf.(3.14)). In 
Section 4 we introduce spatial, temporal and vertical geodesics and state their characterizations via the 
geometric objects defined on ($\bar{M}, \bar{g}$) (cf. Theorems 4.1 and 4.5). In case ${\cal{T}}\bar{M}\oplus 
 {\cal{V}}\bar{M}$ is a Killing vector bundle, we show that spatial geodesics coincide with autoparallel curves
 of $\nabla$ (cf. Theorem 4.3). Finally, we describe explicitly the above three categories of geodesics in a 
$5D$ Robertson-Walker universe (cf. Theorem 4.4 and Corollary 4.2). The conclusions on the research developed 
in the paper are presented in Section 5.
\newpage

\section{Kinematic quantities and the Riemannian spatial connection in a $5D$ universe}

In this section we describe the geometric configuration of a $5D$ universe that has been presented in \cite{b}. 
Let  $\bar{M} = M\times K$ be a product bundle over $M$, where $M$ and $K$ are manifolds of dimensions four and 
one, respectively. The existence of two vector fields $\eta$ and $U$ on $\bar{M}$ and $M$ respectively, induces
 a coordinate system $(x^a)$ on $\bar{M}$ such that $\eta = \frac{\partial }{\partial x^4}$ and $U = 
\frac{\partial }{\partial x^0}$. The coordinate transformations on $\bar{M}$ are given by 

$$\begin{array}{c}
x^{\alpha} = \widetilde{x}^{\alpha}(x^1, x^2, x^3); \ \ \ \widetilde{x}^0 = x^0 + 
f(x^1, x^2, x^3) \vspace{3mm}\\ \ \widetilde{x}^4 = x^4 + \bar{f}(x^0, x^1, x^2, x^3).\end{array}\eqno(2.1)$$

Throughout the paper we use the ranges of indices: $a, b, c, ... \in\{0, 1,2,3,4\},$ $i, j, k, ... \in\{0, 1,2,3\},$ 
$\alpha, \beta, \gamma, ... \in\{1, 2, 3\}$. Also, for any vector bundle $E$ over $\bar{M}$ denote by 
$\Gamma(E)$ the ${\cal{F}}(\bar{M})$-module of smooth sections of $E$, where ${\cal{F}}(\bar{M})$ is the algebra 
of smooth functions on $\bar{M}$.\par
Next, suppose that $\bar{M}$ is endowed with a Lorentz metric $\bar{g}$ such that

$$\bar{g}(\eta, \eta) = \Psi^2. \eqno(2.2) $$
 Denote by ${\cal{V}}\bar{M}$ the line bundle over $\bar{M}$ spanned by $\eta$,  and by ${\cal{H}}\bar{M}$
its complementary orthogonal vector bundle in $T\bar{M}$. Then, suppose that the lift of $\partial/\partial x^0$
 to $\bar{M}$ is timelike with respect to $\bar{g},$ and denote by $\delta/\delta x^0$ its projection
 on ${\cal{H}}(\bar{M})$, that is, we have 

$$\frac{\delta }{\delta x^0} = \frac{\partial }{\partial x^0} - 
A_0\frac{\partial }{\partial x^4}. \eqno(2.3)$$
It is proved that exists a globally defined vector field $\xi$ on $\bar{M}$ which is locally given by 
$\delta/\delta x^0$, and we have

$$\bar{g}(\xi, \xi) = -\Phi^2. \eqno(2.4)$$
Thus the tangent bundle of $\bar{M}$ admits the orthogonal decomposition

$$ T\bar{M} = {\cal{T}}\bar{M}\oplus {\cal{S}}\bar{M}\oplus {\cal{V}}\bar{M},\eqno(2.5)$$
where ${\cal{T}}\bar{M}$ is the line bundle spanned by $\xi$, and ${\cal{S}}\bar{M}$ is the complementary 
orthogonal distribution to ${\cal{T}}\bar{M}$ in ${\cal{H}}\bar{M}$. We call ${\cal{T}}\bar{M}$, 
${\cal{S}}\bar{M}$ and ${\cal{V}}\bar{M}$ the {\it temporal distribution}, the {\it spatial distribution} and 
the {\it vertical distribution}, respectively. According to (2.5) there exists an {\it adapted frame field}
$\{\delta /\delta x^0, \delta /\delta x^{\alpha}, \partial /\partial x^4\}$  on $\bar{M}$, where we put

$$\frac{\delta }{\delta x^{\alpha}} = \frac{\partial }{\partial x^{\alpha}} - B_{\alpha}\frac{\delta }{\delta x^0} 
- A_{\alpha}\frac{\partial }{\partial x^4}. \eqno(2.6)$$ 
Its dual frame field $\{\delta x^0, d x^{\alpha}, \delta x^4 \}$, where we put

$$\delta x^0 = dx^0 + B_{\alpha}dx^{\alpha},  \ \ \ \ \delta x^4 = dx^4 + A_idx^i, \eqno(2.7)$$
is called an {\it adapted coframe field} on $\bar{M}$. The pair $(\bar{M}, \bar{g})$ with the geometric 
configuration presented above is called a $5D$ {\it universe}, and it is the main object of study in the present paper. 

Now, denote by $h$ the Riemannian metric induced by $\bar{g}$ on ${\cal{S}}\bar{M}$, and put 

$$h_{\alpha\beta} = h\left(\frac{\delta}{\delta x^{\beta}}, \frac{\delta}{\delta x^{\alpha}}\right), \ 
\alpha, \beta\in\{1,2,3\}.\eqno(2.8)$$
Then, the line element with respect to the adapted coframe field is given by

$$d\bar{s}^2 = -\Phi^2(\delta x^0)^2 + h_{\alpha\beta}dx^{\alpha}dx^{\beta} + \Psi^2(\delta x^4)^2.\eqno(2.9) $$
The $4D$ {\it vorticity} $\omega_{\alpha\beta}$ and $5D$ {\it vorticity} $\eta_{\alpha\beta}$ in the $5D$ universe 
are given by 

$$\begin{array}{lc} 
\omega_{\alpha\beta} = \frac{1}{2}\left\{\frac{\delta B_{\beta}}{\delta x^{\alpha}} - 
\frac{\delta B_{\alpha}}{\delta x^{\beta}}\right\}, \vspace{2mm}\\ 
\eta_{\alpha\beta} = \frac{1}{2}\left\{\frac{\delta A_{\beta}}{\delta x^{\alpha}} - 
\frac{\delta A_{\alpha}}{\delta x^{\beta}} + B_{\alpha}\frac{\delta A_0}{\delta x^{\beta}} - B_{\beta}
\frac{\delta A_0}{\delta x^{\alpha}}\right\}. \end{array}\eqno(2.10)$$
Also, the $4D$ {\it expansion tensor field} $\Theta_{\alpha\beta}$ and the $5D$ {\it expansion tensor field} 
 $K_{\alpha\beta}$ are given by

$$(a) \ \ \Theta_{\alpha\beta} = \frac{1}{2}\frac{\delta h_{\alpha\beta}}{\delta x^0}, \ \ \ (b) \ \ \ 
K_{\alpha\beta} = \frac{1}{2}\frac{\partial h_{\alpha\beta}}{\partial x^4}. \eqno(2.11)$$
The $4D$ {\it expansion function} $\Theta$ and the $5D$ {\it expansion function} $K$ are the traces of the 
spatial tensor fields from (2.11), expressed as follows 

$$(a) \ \ \ \Theta = \Theta_{\alpha\beta}h^{\alpha\beta}, \ \ \ (b) \ \ \ K 
= K_{\alpha\beta}h^{\alpha\beta}.\eqno(2.12) $$
{\bf Remark 1.1} It is worth mentioning that $\omega_{\alpha\beta}$, $\eta_{\alpha\beta}$, $\Theta_{\alpha\beta}$
and $K_{\alpha\beta}$  define spatial tensor fields of type (0, 2) on the $5D$ universe $(\bar{M}, \bar{g})$. 
According to the study presented in \cite{b}, this means that with respect to the transformations (2.1) they 
change like tensor fields of type (0,2) on a 3-dimensional manifold. $\Box$\par 
 An important geometric object introduced in \cite{b} is the {\it Riemannian spatial connection} on $\bar{M}$, which 
is a linear connection $\nabla$ on the spatial distribution ${\cal{S}}\bar{M}$, given by
 
$$\begin{array}{l}
\nabla_X{\cal{S}}Y = {\cal{S}}\bar{\nabla}_X{\cal{S}}Y, \ \ \ \forall \ \ X,Y 
\in \Gamma(T\bar{M}),\end{array}\eqno(2.13)$$ 
where $\bar{\nabla}$ is the Levi-Civita connection on $\bar{M}$ and ${\cal{S}}$ is the projection morphism
 of $T\bar{M}$ on ${\cal{S}}\bar{M}$ with respect to (2.5). Locally, $\nabla$ is given by

$$\begin{array}{c}
(a) \ \ \nabla_{\frac{\delta }{\delta x^{\beta}}}\frac{\delta }{\delta x^{\alpha}} = 
\Gamma^{\ \gamma}_{{\alpha}\ \;{\beta}}\frac{\delta }{\delta x^{\gamma}}, \ \ \ (b) \ \ \ \ 
\nabla_{\frac{\delta }{\delta x^0}}\frac{\delta }{\delta x^{\alpha}} =  
\Gamma^{\ \gamma}_{\alpha \ \;0}\frac{\delta }{\delta x^{\gamma}}\vspace{2mm}\\ 
(c) \ \ \nabla_{\frac{\partial }{\partial x^4}}\frac{\delta }{\delta x^{\alpha}} = 
\Gamma^{\ \gamma}_{\alpha \ \;4}\frac{\delta }{\delta x^{\gamma}},\end{array}\eqno(2.14) $$
where we put 

$$\begin{array}{lc} (a) \ \ \ 
\Gamma^{\ \gamma}_{\alpha\ \;\beta} = \frac{1}{2}h^{\gamma \mu}\left\{\frac{\delta h_{\mu\alpha}}
{\delta x^{\beta}} +  
\frac{\delta h_{\mu\beta}}{\delta x^{\alpha}} - \frac{\delta h_{\alpha\beta}}{\delta x^{\mu}}\right\},
\vspace{3mm}\\
(b) \ \ \ \Gamma^{\ \gamma}_{\alpha \ \;0} = \Theta_{\alpha}^{\gamma} + \Phi^2\omega_{\alpha}^{\gamma},
\ \ \ (c) \ \ \ \Gamma^{\ \gamma}_{\alpha \ \;4} = K_{\alpha}^{\gamma} 
- \Psi^2\eta_{\alpha}^{\gamma}.\end{array} \eqno(2.15) $$
Next, we express the Lie brackets of vector fields from the adapted frame field, as follows:
 
$$\begin{array}{l}(a) \ \ \ \left [\frac{\delta}{\delta x^{\alpha}}, \frac{\delta}{\delta x^0}\right ] = 
b_{\alpha}\frac{\delta}{\delta x^0} + a_{\alpha}\frac{\partial }{\partial x^4}, \ \ \ (b) \ \ \ 
\left [\frac{\delta}{\delta x^0}, \frac{\partial}{\partial x^4}\right ] = a_0\frac{\partial }{\partial x^4},
\vspace{2mm}\\ (c) \ \ \ \left [\frac{\delta}{\delta x^{\alpha}}, \frac{\partial}{\partial x^4}\right ] = 
d_{\alpha}\frac{\delta}{\delta x^0} + c_{\alpha}\frac{\partial }{\partial x^4},\vspace{2mm}\\ 
 (d) \ \ \ \left [\frac{\delta}{\delta x^{\beta}}, \frac{\delta}{\delta x^{\alpha}}\right ] = 
2\omega_{\alpha\beta}\frac{\delta}{\delta x^0} + 2\eta_{\alpha\beta}\frac{\partial }{\partial x^4}, 
 \end{array}\eqno(2.16)$$
where we put      

$$\begin{array}{lc}a_{\alpha} = \frac{\delta A_{\alpha}}{\delta x^0} - \frac{\delta A_0}{\delta x^{\alpha}}
 -  B_{\alpha}\frac{\delta A_0}{\delta x^0},\vspace{2mm} \\ b_{\alpha} = \frac{\delta B_{\alpha}}{\delta x^0},
 \ \ \ c_{\alpha} = \frac{\partial A_{\alpha}}{\partial x^4} -  B_{\alpha}\frac{\partial A_0}{\partial x^4}, \ \ 
 \ d_{\alpha} = \frac{\partial B_{\alpha}}{\partial x^4}
.\end{array}\eqno(2.17) $$
Note that $a_{\alpha}, b_{\alpha}, c_{\alpha}$ and $d_{\alpha}$ define spatial tensor fields of type $(0,1)$ 
on $\bar{M}$.\par
Finally,  the Levi-Civita connection $\bar{\nabla}$ on $(\bar{M}, \bar{g})$, is expressed as follows:

$$\begin{array}{lc}
\bar{\nabla}_{\frac{\delta }{\delta x^{\beta}}}\frac{\delta }{\delta x^{\alpha}} =  
\Gamma^{\ \;\gamma}_{{\alpha}\ \;{\beta}} \frac{\delta }{\delta x^{\gamma}} + \left(\omega_{{\alpha}{\beta}} 
+ \Phi^{-2}\Theta_{{\alpha}{\beta}}\right)\frac{\delta }{\delta x^0} 
\vspace{2mm}\\ \hspace*{20mm} + \left(\eta_{{\alpha}{\beta}} 
- \Psi^{-2}K_{{\alpha}{\beta}}\right)\frac{\partial }{\partial x^4},\vspace{2mm}\\ 
 
\bar{\nabla}_{\frac{\delta }{\delta x^0}}\frac{\delta }{\delta x^{\alpha}} =  
\Gamma^{\ \;\gamma}_{{\alpha}\ \;0} \frac{\delta }{\delta x^{\gamma}} + \left(\Phi_{{\alpha}} - 
b_{\alpha}\right)\frac{\delta }{\delta x^0} \vspace{2mm}\\ \hspace*{20mm} + 
\frac{1}{2}\left(\Phi^2d_{\alpha}\Psi^{-2} - a_{\alpha}\right)\frac{\partial }{\partial x^4},\vspace{2mm}\\ 

\bar{\nabla}_{\frac{\partial }{\partial x^4}}\frac{\delta }{\delta x^{\alpha}} =  
\Gamma^{\ \;\gamma}_{{\alpha}\ \;4} \frac{\delta }{\delta x^{\gamma}} + 
\frac{1}{2}\left(\Psi^2a_{\alpha}\Phi^{-2} - d_{\alpha}\right)\frac{\delta }{\delta x^0} \vspace{2mm}\\ 
 \hspace*{20mm} + \left(\Psi_{{\alpha}} - c_{\alpha}\right)\frac{\partial }{\partial x^4}, \vspace{2mm}\\ 
 
\bar{\nabla}_\frac{\delta }{\delta x^{\alpha}}{\frac{\delta }{\delta x^0}} =  
\Gamma^{\ \;\gamma}_{{\alpha}\ \;0}\frac{\delta }{\delta x^{\gamma}} + \Phi_{{\alpha}}\frac{\delta }{\delta x^0} + 
\frac{1}{2}\left(\Phi^2d_{\alpha}\Psi^{-2} + a_{\alpha}\right)\frac{\partial }{\partial x^4},\vspace{2mm}\\
 
\bar{\nabla}_{\frac{\delta }{\delta x^{\alpha}}}{\frac{\partial }{\partial x^4}} =  
\Gamma^{\ \;\gamma}_{{\alpha}\ \;4} \frac{\delta }{\delta x^{\gamma}} + 
\frac{1}{2}\left(\Psi^2a_{\alpha}\Phi^{-2} + d_{\alpha}\right)\frac{\delta }{\delta x^0} 
+ \Psi_{{\alpha}}\frac{\partial }{\partial x^4}, \vspace{2mm}\\ 

\bar{\nabla}_{\frac{\partial }{\partial x^4}}\frac{\delta }{\delta x^0} = \frac{1}{2}\left(\Psi^2a^{\gamma} - 
\Phi^{2}d^{\gamma}\right)\frac{\delta }{\delta x^{\gamma}} + \Phi_4\frac{\delta }{\delta x^0} + 
(\Psi_0 - a_0)\frac{\partial }{\partial x^4}, \vspace{2mm}\\  
\bar{\nabla}_\frac{\delta }{\delta x^0}{\frac{\partial }{\partial x^4}} = \frac{1}{2}\left(\Psi^2a^{\gamma} - 
\Phi^{2}d^{\gamma}\right)\frac{\delta }{\delta x^{\gamma}} + \Phi_4\frac{\delta }{\delta x^0} + 
\Psi_0\frac{\partial }{\partial x^4}, \vspace{2mm}\\
\bar{\nabla}_\frac{\delta }{\delta x^0}{\frac{\delta }{\delta x^0}} = \Phi^2\left(\Phi^{\gamma} - 
b^{\gamma}\right)\frac{\delta }{\delta x^{\gamma}} + \Phi_0\frac{\delta }{\delta x^0} + 
\Phi^2\Phi_4\Psi^{-2}\frac{\partial }{\partial x^4}, \vspace{2mm}\\
\bar{\nabla}_\frac{\partial }{\partial x^4}{\frac{\partial }{\partial x^4}} = \Psi^2\left(c^{\gamma} - 
\Psi^{\gamma}\right)\frac{\delta }{\delta x^{\gamma}} + \Psi^2(\Psi_0 - a_0)\Phi^{-2}\frac{\delta }{\delta x^0} + 
\Psi_4\frac{\partial }{\partial x^4}, \end{array}\eqno(2.18)$$
where we put

$$\begin{array}{lc}\Phi_i = \Phi^{-1}\frac{\delta\Phi}{\delta x^i}, \ \ \ 
\Psi_i = \Psi^{-1}\frac{\delta\Psi}{\delta x^i},\vspace{2mm} \\  \Phi_4 = \Phi^{-1}\frac{\partial\Phi}{\partial x^4},  
\ \ \ \Psi_4 = \Psi^{-1}\frac{\partial\Psi}{\partial x^4}.\end{array}\eqno(2.19)$$

\section{Equations of motion in a $5D$ universe}

In this section we write down, in a covariant form, the equations of motion in the $5D$ universe $(\bar{M}, 
\bar{g})$. It is first time in literature when these equations are expressed by three groups of equations (cf. 
(3.6)), and in terms of kinematic quantities and of the local coefficients of the Riemannian spatial connection. 
As an example of such $5D$ universe we present the $5D$ Robertson-Walker universe (cf.(3.7)), and state its
 equations of motion (cf.(3.14)).\par 
Let $\bar{C}$ be a smooth curve in $\bar{M}$ given by the equations 

$$x^a = x^a(t), \ \ \ t\in [c,d], \ \ \ a\in\{0,1,2,3,4\}.\eqno(3.1)$$
Then by direct calculations using (2.3) and (2.6), we deduce that the tangent vector field $\frac{d }{dt}$ to 
$\bar{C}$ is expressed with respect to the adapted frame field $\{\frac{\delta}{\delta x^0}, \frac{\delta}{\delta x^
{\alpha}}, \frac{\partial }{\partial x^4}\}$, as follows
 
$$\frac{d}{dt} = \frac{\delta x^0}{\delta t}\frac{\delta }{\delta x^0} + \frac{dx^{\alpha}}{dt}\frac{\delta }
{\delta x^{\alpha}} + \frac{\delta x^4}{\delta t}\frac{\partial }{\partial  x^4}, \eqno(3.2)$$ 
where we put 

$$\frac{\delta x^0}{\delta t} = \frac{d x^0}{d t} + B_{\alpha}\frac{d x^{\alpha}}{d t}, \ \ \ 
\frac{\delta x^4}{\delta t} = \frac{d x^4}{d t} + A_i\frac{d x^i}{d t}.\eqno(3.3)$$
Next, after some long calculations by using (3.2) and (2.18) we obtain

$$\begin{array}{lc}
\bar{\nabla}_{\frac{d}{dt}}\frac{\delta }{\delta x^0} = \left\{\Phi^2\left (\Phi^{\gamma} 
- b^{\gamma}\right )\frac{\delta x^0}{\delta t}\right.\vspace{2mm}\\ 
 \hspace*{16mm} + \left. \Gamma^{\;\;\gamma}_{\alpha \ 0}\frac{dx^{\alpha}}{dt} + 
\frac{1}{2}\left ( \Psi^2a^{\gamma} - \Phi^2d^{\gamma}\right )\frac{\delta x^4}{\delta t}\right\}\frac{\delta }
{\delta x^{\gamma}}\vspace{2mm}\\  \hspace*{16mm} + \left\{\Phi_0\frac{\delta x^0}{\delta t} + 
 \Phi_{\alpha}\frac{dx^{\alpha}}{dt} + \Phi_4\frac{\delta x^4}{dt} \right\}\frac{\delta}{\delta x^0} 
\vspace{2mm}\\ \hspace*{16mm} + \left\{\Phi^2\Phi_4\Psi^{-2}\frac{\delta x^0}{\delta t} + 
\frac{1}{2}\left (\Phi^2d_{\alpha}\Psi^{-2} + a_{\alpha}\right )\frac{dx^{\alpha}}{dt} \right. \vspace{2mm}\\ 
 \hspace*{16mm} \left. + \left (\Psi_0 - a_0\right )\frac{\delta x^4}{\delta t} \right\}
\frac{\partial }{\partial x^4},\vspace{2mm}\\
 
 \bar{\nabla}_{\frac{d}{dt}}\frac{\delta }{\delta x^{\alpha}} = \left\{\Gamma^{\;\;\gamma}_{\alpha \ 0}
\frac{\delta x^0}{\delta t} +  \Gamma^{\;\;\gamma}_{\alpha \ \beta}\frac{dx^{\beta}}{dt} + \Gamma^{\;\;\gamma}_
{\alpha \ 4}\frac{\delta x^4}{\delta t}\right\}\frac{\delta }{\delta x^{\gamma}} \vspace{2mm}\\ \hspace*{16mm} + 
\left\{\left (\Phi_{\alpha} - b_{\alpha}\right )\frac{\delta x^0}{\delta t} + \left (\omega_{\alpha\beta} + 
\Phi^{-2}\Theta_{\alpha\beta}\right )\frac{dx^{\beta}}{dt}\right.\vspace{2mm}\\ \hspace*{16mm} + \left. 
 \frac{1}{2}\left (\Psi^2a_{\alpha}\Phi^{-2} - d_{\alpha} \right )\frac{\delta x^4}{\delta t}\right\}
\frac{\delta }{\delta x^0}\vspace{2mm}\\ \hspace*{16mm} + \left\{\frac{1}{2}\left (\Phi^2d_{\alpha}\Psi^{-2}
 - a_{\alpha} \right )\frac{\delta x^0}{\delta t} + \left (\eta_{\alpha\beta} - \Psi^{-2}K_{\alpha\beta}
\right )\frac{dx^{\beta}}{dt} \right.\vspace{2mm}\\ \hspace*{16mm} + \left. 
\left (\Psi_{\alpha} - c_{\alpha}\right )\frac{\delta x^4}{\delta t} \right\}\frac{\partial }
{\partial x^4},\vspace{2mm}\\

\bar{\nabla}_{\frac{d}{dt}}\frac{\partial }{\partial x^4} = \left\{\frac{1}{2}\left (\Psi^2a^{\gamma} - 
\Phi^2d^{\gamma}\right )\frac{\delta x^0}{\delta t} +  \Gamma^{\;\;\gamma}_{\alpha \ 4}\frac{dx^{\alpha}}{dt}
 \right. \vspace{2mm}\\ \hspace*{14mm} \left. + \Psi^2\left ( c^{\gamma} - \Psi^{\gamma}\right )
\frac{\delta x^4}{\delta t}\right\}\frac{\delta }{\delta x^{\gamma}} + \left\{\Phi_4\frac{\delta x^0}{\delta t}
\right. \vspace{2mm}\\ \hspace*{16mm} \left. + \frac{1}{2}\left ( \Psi^2a_{\alpha} \Phi^{-2} + d_{\alpha}
\right) \frac{dx^{\alpha}}{dt} + \Psi^2\left (\Psi_0 - a_0\right )\frac{\delta x^4}{\delta t} \right\}
\frac{\delta }{\delta x^0}\vspace{2mm}\\ \hspace*{16mm} + \left\{\Psi_0\frac{\delta x^0}{\delta t}
 + \Psi_{\alpha}\frac{dx^{\alpha}}{dt} + \Psi_4\frac{\delta x^4}{\delta t}\right\}\frac{\partial }{\partial x^4}.
\end{array}\eqno(3.4)$$
Now, by using (3.2), (3.4), (2.15b) and (2.15c), and taking into account that $\omega_{\alpha\beta}$ and 
 $\eta_{\alpha\beta}$ are skew-symmetric spatial tensor fields on $\bar{M}$, we deduce that

$$\begin{array}{lc}
\bar{\nabla}_{\frac{d}{dt}}\frac{d }{dt} = \frac{d^2x^{\gamma}}{dt^2}\frac{\delta }{\delta x^{\gamma}} + \frac{d}{dt}
\left(\frac{\delta x^0}{\delta t}\right )\frac{\delta }{\delta x^0} + \frac{d}{dt}
\left(\frac{\delta x^4}{\delta t}\right )\frac{\partial }{\partial x^4} + \frac{dx^{\alpha}}{dt}
\bar{\nabla}_{\frac{d}{dt}}\frac{\delta }{\delta x^{\alpha}}\vspace{2mm}\\ \hspace*{12mm}

 +  \frac{\delta x^0}{\delta t}
\bar{\nabla}_{\frac{d}{dt}}\frac{\delta }{\delta x^0} + \frac{\delta x^4}{\delta t}\bar{\nabla}_{\frac{d}{dt}}
\frac{\partial}{\partial x^4} = \left\{\frac{d^2x^{\gamma}}{dt^2} + \Gamma^{\ \gamma}_{\alpha 
 \;\;\beta}\frac{dx^{\alpha}}{dt}\frac{dx^{\beta}}{dt}\right. \vspace{2mm}\\ \hspace*{12mm} \left. + 
 
2\left (\Theta^{\gamma}_{\alpha} + \Phi^2\omega_{\alpha}^{\gamma}\right )\frac{dx^{\alpha}}{dt}
\frac{\delta x^0}{\delta t} + 2\left (K_{\alpha}^{\;\;\gamma} - \Psi^2\eta_{\alpha}^{\;\gamma}\right ) 
 \frac{dx^{\alpha}}{dt}\frac{\delta x^4}{\delta t}\right. \vspace{2mm}\\ \hspace*{12mm} \left.  + 
 
\left (\Psi^2a^{\gamma} - \Phi^2d^{\gamma}  \right )\frac{\delta x^0}{\delta t}\frac{\delta x^4}{\delta t}
 + \Phi^2\left (\Phi^{\gamma} - b^{\gamma}\right ) \left (\frac{\delta x^0}{\delta t}\right )^2\right. 
\vspace{2mm}\\ \hspace*{12mm} \left. 

+ \Psi^2\left ( c^{\gamma} - \Psi^{\gamma}\right ) \left (
\frac{\delta x^4}{\delta t}\right )^2 \right\}\frac{\delta }{\delta x^{\gamma}} + \left\{\frac{d}{dt}
\left ( \frac{\delta x^0}{\delta t}\right ) + 
 \Phi^{-2}\Theta_{\alpha\beta}\frac{dx^{\alpha}}{dt}\frac{dx^{\beta}}{dt}\right. 
\vspace{2mm}\\ \hspace*{12mm} \left. + 

\left (2\Phi_{\alpha} - b_{\alpha}\right )\frac{dx^{\alpha}}{dt}
\frac{\delta x^0}{dt} + \Psi^2a_{\alpha}\Phi^{-2}\frac{dx^{\alpha}}{dt}\frac{\delta x^4}{dt} + 2\Phi_4
\frac{\delta x^0}{\delta t}\frac{\delta x^4}{dt}\right. \vspace{2mm}\\ \hspace*{12mm} \left. + \Phi_0
\left (\frac{\delta x^0}{\delta t}\right )^2 + \Psi^2\left (\Psi_0 - a_0 \right )
\left (\frac{\delta x^4}{\delta t}\right )^2\right\}\frac{\delta }{\delta x^0} + \left\{\frac{d}{dt}
\left (\frac{\delta x^4}{\delta t}\right )\right. \vspace{2mm}\\ \hspace*{12mm} \left. 

- \Psi^{-2}K_{\alpha\beta}\frac{dx^{\alpha}}{dt}\frac{dx^{\beta}}{dt}
 + \Phi^2d_{\alpha}\Psi^{-2}\frac{dx^{\alpha}}{dt}\frac{\delta x^0}{\delta t} + \left (2\Psi_{\alpha} - 
c_{\alpha} \right )\frac{dx^{\alpha}}{dt}\frac{\delta x^4}{\delta t}\right. \vspace{2mm}\\ \hspace*{12mm} 
\left.  + 

\left (2\Psi_0 - a_0 \right )\frac{\delta x^0}{\delta t}\frac{\delta x^4}{\delta t}
 + \Phi^2\Phi_4\Psi^{-2}\left (\frac{\delta x^0}{\delta t}\right )^2 + \Psi_4\left (\frac{\delta x^4}{\delta t}
\right )^2\right\}\frac{\partial }{\partial x^4}.\end{array}\eqno(3.5)$$
Finally, since $\bar{C}$ is a geodesic of $(\bar{M}, \bar{g})$, if and only if, the left hand side in (3.5) 
vanishes identically on $\bar{M}$, we can state the following theorem.\vspace{4mm}\par

{\bf Theorem 3.1} {\it Let $(\bar{M}, \bar{g})$ be a $5D$ universe with kinematic quantities $\{\omega_{\alpha\beta},
 \eta_{\alpha\beta}, \Theta_{\alpha\beta}, K_{\alpha\beta}\}$ and with the Riemannian spatial connection $\nabla$ 
 given by} (2.14) {\it and} (2.15). {\it Then the equations of motion in  $(\bar{M}, \bar{g})$ are expressed
 as follows}:

$$\begin{array}{lc}(a) \ \ \ 
\frac{d^2x^{\gamma}}{dt^2} + \Gamma^{\ \gamma}_{\alpha 
 \;\;\beta}\frac{dx^{\alpha}}{dt}\frac{dx^{\beta}}{dt} + 
2\left (\Theta^{\gamma}_{\alpha} + \Phi^2\omega_{\alpha}^{\gamma}\right )\frac{dx^{\alpha}}{dt}
\frac{\delta x^0}{\delta t} + 2\left (K_{\alpha}^{\;\;\gamma}\right.\vspace{2mm}\\ \hspace*{6mm} \left. 

- \Psi^2\eta_{\alpha}^{\;\gamma}\right ) \frac{dx^{\alpha}}{dt}\frac{\delta x^4}{\delta t}  + 
\left (\Psi^2a^{\gamma} - \Phi^2d^{\gamma}  \right )\frac{\delta x^0}{\delta t}\frac{\delta x^4}{\delta t} 
 + \Phi^2\left (\Phi^{\gamma} \right.\vspace{2mm}\\ \hspace*{6mm} \left. 

- b^{\gamma}\right ) \left (
\frac{\delta x^0}{\delta t}\right )^2   + \Psi^2\left ( c^{\gamma} - \Psi^{\gamma}\right ) \left (
\frac{\delta x^4}{\delta t}\right )^2 = 0, \ \ \mbox{where}\ \ \gamma\in\{1,2,3\},\vspace{2mm}\\ 
(b) \ \ \  \frac{d}{dt}\left ( \frac{\delta x^0}{\delta t}\right ) + 
 \Phi^{-2}\Theta_{\alpha\beta}\frac{dx^{\alpha}}{dt}\frac{dx^{\beta}}{dt} 
+ \left (2\Phi_{\alpha} - b_{\alpha}\right )\frac{dx^{\alpha}}{dt}\frac{\delta x^0}{dt}\vspace{2mm}\\ 
\hspace*{12mm} + \Psi^2a_{\alpha}\Phi^{-2}\frac{dx^{\alpha}}{dt}\frac{\delta x^4}{dt} + 2\Phi_4
\frac{\delta x^0}{\delta t}\frac{\delta x^4}{dt} + \Phi_0\left (\frac{\delta x^0}{\delta t}\right )^2
\vspace{2mm}\\ \hspace*{12mm} + \Psi^2\left (\Psi_0 - a_0 \right )\left (\frac{\delta x^4}{\delta t}\right )^2
 = 0,\vspace{2mm}\\ 
(c) \ \ \ \frac{d}{dt}
\left (\frac{\delta x^4}{\delta t}\right ) - \Psi^{-2}K_{\alpha\beta}\frac{dx^{\alpha}}{dt}\frac{dx^{\beta}}{dt}
 + \Phi^2d_{\alpha}\Psi^{-2}\frac{dx^{\alpha}}{dt}\frac{\delta x^0}{\delta t}\vspace{2mm}\\ \hspace*{6mm} + 
\left (2\Psi_{\alpha} - c_{\alpha} \right )\frac{dx^{\alpha}}{dt}\frac{\delta x^4}{\delta t} + \left (2\Psi_0 
- a_0 \right )\frac{\delta x^0}{\delta t}\frac{\delta x^4}{\delta t} + \Phi^2\Phi_4\Psi^{-2}\left (
\frac{\delta x^0}{\delta t}\right )^2\vspace{2mm}\\ \hspace*{6mm} + \Psi_4\left (\frac{\delta x^4}{\delta t}
\right )^2 = 0.\end{array}\eqno(3.6)$$
It is the first time in literature when the equations of motion in a $5D$ universe are expressed in terms of 
kinematic quantities and of some spatial tensor fields. The first type of equations of motion was presented in 
formula (5.28) of  \cite{w}, wherein the natural frame field $\{\partial/\partial x^a\}, \ \ a\in\{0,1,2,3,4 \}$,
 has been used. In this way, no differences were noticed between temporal variable $x^0$, the spatial variable 
 $(x^{\alpha})$ and the vertical variable $x^4$. Also, in \cite{ab} the author stated another form of equations
 of motion, wherein the temporal distribution was not taken into consideration. The main difference between (5.6)
 of \cite{ab} and (3.6), is that the latter can relate physics and geometry with observations, via the kinematic 
quantities.\par
Next, we construct an example of $5D$ universe and write down its equations of motion. Suppose that the line
 element of the Lorentz metric $\bar{g}$ has the particular form

$$d\bar{s}^2 = -(dx^0)^2 + f^2(x^0, x^4)g_{\alpha\beta}(x^1, x^2, x^3)dx^{\alpha}dx^{\beta} + (dx^4)^2,
\eqno(3.7)$$
where $f$ is a positive smooth function on an open region ${\cal{R}}$ of $R_1^2$, and $g_{\alpha\beta}$ define
 a positive definite symmetric spatial tensor field $g$ on $\bar{M}$. Taking into account thst $\bar{g}$ given
 by (3.7) satisfies 

$$\bar{g}(\frac{\partial }{\partial x^0}, \frac{\partial }{\partial x^4}) = 0, \ \ \ \bar{g}(\frac{\partial }
{\partial x^{\alpha}}, \frac{\partial }{\partial x^4}) = 0, \ \ \ \bar{g}(\frac{\partial }{\partial 
x^{\alpha}}, \frac{\partial }{\partial x^0}) = 0, $$
and using (2.3) and (2.6), we obtain

$$\begin{array}{lc}
(a) \ \ \frac{\delta }{\delta x^i} = \frac{\partial }{\partial x^i}, \ \ (b) \ \ \ A_i = 0, \vspace{2mm} \\ (c) \ \ 
B_{\alpha} = 0, \ \forall i\in\{0, 1, 2, 3\},\ \ \alpha\in\{1,2,3\}.\end{array}\eqno(3.8)
$$
By using (3.8a) we see that the distributions ${\cal{S}}\bar{M}$, ${\cal{T}}\bar{M}\oplus {\cal{S}}\bar{M}$ 
and ${\cal{S}}\bar{M}\oplus {\cal{V}}\bar{M}$ are integrable, and as a consequence of (2.16) we deduce that 

$$\begin{array}{lc}
a_i = 0, \ \ b_{\alpha} = c_{\alpha} = d_{\alpha} = 0, \\
\omega_{\alpha\beta} = \eta_{\alpha\beta} = 0, \ \forall i\in\{0, 1, 2, 3\},\ \ \alpha, \beta\in\{1,2,3\}.
\end{array} \eqno(3.9)$$
In order to obtain the other kinematic quantities, we note that 

$$h_{\alpha\beta} = f^2g_{\alpha\beta} \ \ \mbox{and} \ \ h^{\alpha\beta} = f^{-2}g^{\alpha\beta}.\eqno(3.10)
$$
Then, by using (3.10) into (2.11) and (2.12), we infer that 

$$\begin{array}{lc}\Theta_{\alpha\beta} = f\frac{\partial f}{\partial x^0}g_{\alpha\beta}, \ \ \ K_{\alpha\beta}
 = f\frac{\partial f}{\partial x^4}g_{\alpha\beta},\vspace{2mm} \\ \Theta_{\alpha}^{\gamma} = f^{-1}
\frac{\partial f}{\partial x^0}\delta_{\alpha}^{\gamma}, \ \ \ K_{\alpha}^{\gamma} = f^{-1}
\frac{\partial f}{\partial x^4}\delta_{\alpha}^{\gamma},\vspace{2mm} \\ \Theta = 3f^{-1}
\frac{\partial f}{\partial x^0}, \ \ \ K = 3f^{-1}\frac{\partial f}{\partial x^4}.
\end{array}\eqno(3.11)$$
Moreover, from (2.19) and (3.3), we obtain 

$$\Phi_a = \Psi_a = 0, \ \ \ \forall a\in \{0,1,2,3,4\},\eqno(3.12)$$
and

$$\frac{\delta x^0}{\delta t} = \frac{dx^0}{dt}, \ \ \ \frac{\delta x^4}{\delta t} = \frac{dx^4}{dt},\eqno(3.13)$$
respectively. Also, note that the local coefficients $\Gamma_{\alpha\ \ \beta}^{\ \gamma}$ of the Riemannian 
spatial connection $\nabla$ given by (2.15a) become the usual Christoffel symbols with respect to the Riemannian 
 metric $g = (g_{\alpha\beta})$.\newline
Finally, by using (3.7) - (3.13) into (3.6), we obtain the following equations of motion in a $5D$ universe 
$(\bar{M}, \bar{g})$ whose Lorentz metric is given by (3.7):
$$\begin{array}{lc}(a) \ \ \ 
\frac{d^2x^{\gamma}}{dt^2} + \Gamma^{\ \gamma}_{\alpha 
 \;\;\beta}\frac{dx^{\alpha}}{dt}\frac{dx^{\beta}}{dt} + 
2f^{-1}\frac{df}{dt}\frac{d x^{\gamma}}{d t} = 0, \vspace{2mm}\\ 
(b) \ \ \  \frac{d^2x^0}{d t^2} +  f\frac{\partial f}{\partial x^0}g_{\alpha\beta}\frac{dx^{\alpha}}{dt}
\frac{d x^{\beta}}{dt} = 0,\vspace{2mm}\\ 
(c) \ \ \ \frac{d^2 x^4}{d t^2} - f\frac{\partial f}{\partial x^4}g_{\alpha\beta}\frac{dx^{\alpha}}{dt}
\frac{d x^{\beta}}{dt} = 0.\end{array}\eqno(3.14)$$
As leaves of ${\cal{T}}\bar{M}\oplus {\cal{S}}\bar{M}$ are locally given by $x^4 = $const., from (3.7) we see 
that the metric induced on them is of Robertson-Walker metric type (cf.\cite{bo}, p.343), provided the leaves
 of ${\cal{S}}\bar{M}$ are 3-dimensional manifolds of the same constant curvature. This happens in the case we take 
$\bar{M} = I\times S\times K$, where $I$ is an open interval in $R$ and $S$ is a 3-dimensional Riemannian manifold
 of constant curvature $k = 1, 0$ or $-1$. Thus, we may think such a $5D$ universe as a disjoint union of 
 Robertson-Walker spacetimes. For this reason we call the $5D$ universe $(\bar{M}, \bar{g})$ whose metric is given
 by $(3.7)$, a {\it $5D$ Robertson-Walker universe}, with the {\it warping function f}.

\section{Special geodesics in a $5D$-universe}

This section is devoted to the study of some particular classes of geodesics in a  $5D$ universe $(\bar{M}, \bar{g})$. 
The existence of these geodesics is due to the splitting (2.5) of $T\bar{M}$, which has been considered first in 
\cite{b}.\par 
Let $\bar{C}$ be a curve in $\bar{M}$  given by (3.1). Then, we say that $\bar{C}$ is a {\it spatial curve}, if 
it is tangent to the spatial distribution at any of its points. Thus, by (3.2) and (3.3), we deduce that $\bar{C}$
 is a spatial curve, if and only if, we have 

$$\frac{d }{dt} = \frac{dx^{\alpha}}{dt}\frac{\delta }{\delta x^{\alpha}}, \eqno(4.1)$$
or equivalently

$$\frac{\delta x^0}{\delta t} = \frac{dx^0}{dt} + B_{\alpha}\frac{dx^{\alpha}}{dt} = 0, \ \ \ \mbox{and} \ \ \ 
\frac{\delta x^4}{\delta t} = \frac{dx^4}{dt} + A_i\frac{dx^i}{dt} = 0.\eqno(4.2)$$
If moreover, a spatial curve $\bar{C}$ is a geodesic of $(\bar{M}, \bar{g})$, we say that it is a {\it 
spatial geodesic}. Taking into account of (4.2) into (3.6), we state the following theorem.\vspace{4mm}\par
{\bf Theorem 4.1} {\it A spatial curve $\bar{C}$ is a spatial geodesic, if and only if, the following equations
 are satisfied}:

$$\begin{array}{lc} (a) \ \ \ \frac{d^2x^{\gamma}}{dt^2} + \Gamma^{\ \gamma}_{\alpha \;\;\beta} 
 \frac{dx^{\alpha}}{dt}\frac{dx^{\beta}}{dt} = 0,\vspace{2mm}\\ 
(b) \ \ \ \Theta_{\alpha\beta} \frac{dx^{\alpha}}{dt}\frac{dx^{\beta}}{dt} = 0,\vspace{2mm}\\
(c) \ \ \ K_{\alpha\beta} \frac{dx^{\alpha}}{dt}\frac{dx^{\beta}}{dt} = 0.
\end{array}\eqno(4.3)$$
Next, we say that $\bar{C}$ is an {\it autoparallel curve} in $\bar{M}$ with respect to the Riemannian 
spatial connection $\nabla$, if it is a spatial curve satisfying 

$$\nabla_{\frac{d}{dt}}\frac{d}{dt} = 0.\eqno(4.4)$$
Then, by using (4.1) and (2.14a) into (4.4) we obtain the following.\vspace{4mm}\par
{\bf Theorem 4.2} {\it A spatial curve $\bar{C}$ is an autoparallel curve with respect to $\nabla$, if and only
 if, the equations} (4.3a) {\it are satisfied}.\par 
Thus, the relationship between spatial geodesics and autoparallel curves with respect to $\nabla$, can be stated
 in the next corollary.\vspace{4mm}\par
 
{\bf Corollary 4.1} {\it A spatial geodesic of $(\bar{M}, \bar{g})$ must be an autoparallel curve with respect
 to $\nabla$. Conversely, an autoparallel curve with respect to $\nabla$ is a spatial geodesic, if and only if,}
 (4.3b) {\it and} (4.3c) {\it are satisfied}.\par
Next, we define  the Lie derivative of the Riemannian metric $h$ on ${\cal{S}}\bar{M}$ with respect to a vector
 field $Z\in \Gamma({\cal{T}}\bar{M}\oplus{\cal{V}}\bar{M})$ as follows:

$$\begin{array}{lc}
({\cal{L}}_Zh)({\cal{S}}X, {\cal{S}}Y) = Z(h({\cal{S}}X, {\cal{S}}Y)) - h({\cal{S}}[Z, 
{\cal{S}}X], {\cal{S}}Y)\vspace{2mm}\\ \hspace*{32mm} - h({\cal{S}}[Z, {\cal{S}}Y], {\cal{S}}X), 
\forall X, Y \in\Gamma(T\bar{M}).\end{array}\eqno(4.5)$$
Then, take in turn $Z = \delta/\delta x^0$ and $Z = \partial/\partial x^4$ in (4.5), and by using (2.8), 
(2.11), (2.16a) and (2.16c), we obtain 

$$\begin{array}{lc}(a) \ \ \ ({\cal{L}}_{\frac{\delta}{\delta x^0}}h)(\frac{\delta}{\delta x^{\alpha}}, 
 \frac{\delta}{\delta x^{\beta}}) = 2\Theta_{\alpha\beta}, \vspace{2mm}\\ 
(b) \ \ \ ({\cal{L}}_{\frac{\delta}{\delta x^4}}h)(\frac{\delta}{\delta x^{\alpha}}, 
\frac{\delta}{\delta x^{\beta}}) = 2K_{\alpha\beta}.
\end{array}\eqno(4.6)$$
Now, we say that  ${\cal{T}}\bar{M}\oplus {\cal{V}}\bar{M}$ is a {\it Killing vector bundle} with respect to 
$( {\cal{S}}\bar{M}, h)$, if the Lie derivative given by (4.5) vanishes identically on $\bar{M}$, for any $Z$. 
Then, from (4.6) we see that  ${\cal{T}}\bar{M}\oplus {\cal{V}}\bar{M}$ {\it is a Killing vector bundle, if and
 only if, both the $4D$ and $5D$ expansion tensor fields vanish identically on} $\bar{M}.$ Thus, combining 
Theorems 4.1 and 4.2, we obtain the following theorem.\vspace{4mm}\par

{\bf Theorem 4.3} {\it Let $(\bar{M}, \bar{g})$ be a $5D$ universe such that ${\cal{T}}\bar{M}\oplus {\cal{V}}
\bar{M}$ is a Killing vector bundle. Then a spatial curve $\bar{C}$ in $\bar{M}$ is a spatial geodesic, if and
 only if, it is an autoparallel with respect to the Riemannian spatial connection} $\nabla$. \par

In particular, consider  a $5D$ Robertson-Walker universe, and by using (3.8), (3.14) and (4.2), we obtain
 the following.\vspace{4mm}\par

{\bf Theorem 4.4} {\it Let $(\bar{M}, \bar{g})$ be a $5D$ Robertson-Walker universe whose metric is given by}
 (3.7). {\it Then a curve $\bar{C}$ in  $\bar{M}$ is a spatial geodesic, if and only if, the following conditions
 are satisfied}:\par
 (i) \ {\it The parametric equations of $\bar{C}$ have the form}

$$x^0 = c, \ \ \ x^{\gamma} = x^{\gamma}(t), \ \ \ x^4 = k,\eqno(4.7) $$
{\it where $c$ and $k$ are constants, and $x^{\gamma} = x^{\gamma}(t), \ \gamma\in\{1,2,3\}$, define a geodesic 
of a leaf of ${\cal{S}}\bar{M}$ with respect to the Riemannian metric} $g = (g_{\alpha\beta})$.\par
(ii) \ {\it The warping function $f$ admits $(c,k)$ as critical point, that is,}

$$\frac{\partial f}{\partial x^0}(c, k) =  \frac{\partial f}{\partial x^4}(c, k) = 0.\eqno(4.8)$$
The above theorem says that spatial geodesics in $(\bar{M}, \bar{g})$ exist, if and only if, the warping function
 has at least one critical point $(c, k)$. In that case, if $S$ is the leaf of ${\cal{S}}\bar{M}$ given by 
equations $x^0 = c, x^4 = k$, then the lifts of geodesics of $(S, g)$ are spatial geodesics of $(\bar{M}, \bar{g})$.
\par
Finally, we say that a geodesic $\bar{C}$ of $(\bar{M}, \bar{g})$ is a {\it temporal geodesic} (resp. {\it 
vertical geodesic}) if it is tangent to ${\cal{T}}\bar{M}$ (resp. ${\cal{V}}\bar{M}$) at any of its points. Then,
 by using (3.2), (3.3) and (3.6) we state the following theorem.\vspace{4mm}\par

{\bf Theorem 4.5} (i) \ \ {\it A curve $\bar{C}$ is a temporal geodesic in the $5D$ universe $(\bar{M}, \bar{g})$, 
if and only if, the following conditions are satisfied}:

$$\begin{array}{lc}
(a) \ \ \ \Phi_4 = 0, \ \ \ (b) \ \ \ \Phi_{\alpha} = b_{\alpha}, \ \ \ (c) \ \ \ x^{\alpha} = k^{\alpha}, \ 
\ \alpha\in\{1,2,3\}, \\ 
(d) \ \ \ \frac{d^2x^0}{dt^2} + \Phi_0\left(\frac{dx^0}{dt}\right)^2 = 0, \ \ \ (e) \ \ \ \frac{dx^4}{dt} 
+ A_0\frac{dx^0}{dt} = 0,\end{array} \eqno(4.9)$$
{\it where $k^{\alpha}$ are constants.}\par
(ii) \ {\it A curve $\bar{C}$ is a vertical geodesic in  $(\bar{M}, \bar{g})$, 
if and only if, we have}:\par

$$\begin{array}{lc}
(a) \ \ \ \Psi_0 = a_0, \ \ \ (b) \ \ \ \Psi_{\alpha} = c_{\alpha}, \ \\ (c) \ \ \ x^i = \lambda^i, \ 
\ \alpha\in\{1,2,3\}, \ \ i\in\{0,1,2,3\},\\ 
(d) \ \ \ \frac{d^2x^4}{dt^2} + \Psi_4\left(\frac{dx^4}{dt}\right)^2 = 0, \end{array} \eqno(4.10)$$
{\it where $\lambda^i$ are constants.}\par
In general, (4.9a) and (4.9b) (resp. (4.10a) and (4.10b)) are strong constraints which should be satisfied by
 the solutions from (4.9c), (4.9d) and (4.9e) (resp. (4.10c) and (4.10d)). However, we note that all these 
constraints are satisfied in case of a $5D$ Robertson-Walker universe. More precisely, from Theorem 4.5 we deduce
 the following corollary.\vspace{4mm}\par

{\bf Corollary 4.2} {\it Let $(\bar{M}, \bar{g})$ be a $5D$ Robertson-Walker universe. Then we have the 
following assertions}:\par
(i) {\it The temporal geodesics of $(\bar{M}, \bar{g})$ exists, and they are portions of lines given by $x^u = k^u, \ 
u\in\{1,2,3,4\}$, where $k^u$ are constants.}\par
(ii) {\it The vertical geodesics of $(\bar{M}, \bar{g})$ exists, and they are portions of lines given by $x^i =
\lambda^i, \ i\in\{0,1,2,3\}$, where $\lambda^i$ are constants.}

\section{Conclusions}

The present paper has its roots in \cite{b} and \cite{bc}, wherein we developed new approaches on the 
$(1 + 1 + 3)$ threading of a $5D$ universe and on the $(1 + 3)$ threading of a spacetime, respectively. The 
main geometric objects used in the paper are: the adapted frame and coframe fields, the kinematic tensor fields,
 and the Riemannian spatial connection. By using these geometric objects, we state in a $5D$ covariant form, the 
equations of motion in  $(\bar{M}, \bar{g})$. The splitting of such equations in three groups (see (3.6)) enables
 us to consider the spatial, temporal and vertical geodesics. We note the interrelations between spatial geodesics
 and autoparallel curves with respect to the Riemannian spatial connection (cf. Corollary 4.1). In particular, if 
${\cal{T}}\bar{M}\oplus {\cal{V}}\bar{M}$ is a Killing vector bundle, we show that spatial geodesics coincide 
with autoparallel curves of $\nabla$ (cf. Theorem 4.3). This shows that $\nabla$ has an important role in the study
 of geometry and physics of a $5D$ universe.\par 
As a new example of $5D$ universe in the sense considered in \cite{b},
 we present what we call the $5D$ Robertson-Walker universe, whose metric is given by (3.7). We show that such a 
 universe can be thought as a disjoint union of $4D$ Robertson-Walker spacetimes. The equations of motion have 
the simple form (cf. (3.14)), wherein the first two groups remind us of the equations of motion in a $4D$ 
 Robertson-Walker spacetime (cf.\cite{bo}, p.353). Also, we show that the projections of spatial geodesics of 
 $(\bar{M}, \bar{g})$ on the leaves of ${\cal{S}}\bar{M}$ are just geodesics of the leaves with the Riemannian 
metric $g$ (cf. Theorem 4.4).\par
Finally, we note that throughout the paper, the spatial tensor fields enable us to apply the principle of 
covariance, which is one of the most powerful ideas in modern physics, This will be seen more evidently in
 a forthcoming paper on the splitting of Einstein equations in a $5D$ universe.


\begin{thebibliography}{nn}

\bibitem{b} A.Bejancu, Eur.Phys. J.C(2015){\bf 75}:346
\bibitem{w} P.S.Wesson, {\it Space-Time-Matter. Modern Kaluza-Klein Theory} (World Sicentific, Singapore, 1999)
\bibitem{o} J. M.Overduin, P.S.Wesson, Phys. Rep. {\bf 283}, 303(1999) 
\bibitem{m} R.Maartens, K.Koyama, Living Rev.Relativ.{\bf 13}, 5(2010)
\bibitem{l} D.Langlois, Prog.Theor.Phys.Suppl.{\bf 148}, 181(2002)
\bibitem{ab} A.Bejancu, Gen. Relativ.Gravit. {\bf 45}, 2273(2013)
\bibitem{bo} B.O'Neill, {\it Semi- Riemannian Geometry and Applications to Relativity} (Academic Press,
             New York, 1983)
\bibitem{bc} A.Bejancu,  C.C\u alin, Eur.Phys.J.C( 2015 ) {\bf 75}:159
 
\end{thebibliography}
\end{document}